\documentclass[11pt]{article}
\usepackage{amssymb,amsmath}
\usepackage{graphicx,units,mathrsfs,bm}
\usepackage{subfigure}
\usepackage[sort&compress,numbers]{natbib} 
\usepackage{hyperref,xcolor}
\usepackage[paperwidth=8.5in, paperheight=11in, portrait, top=1in, bottom=1.in, left=1.15in, right=1.15in]{geometry}
\usepackage{tikz}
\usetikzlibrary{calc,decorations.markings,positioning}
\allowdisplaybreaks
\newcommand\bt{\raise 2pt \hbox{$\bigtriangledown$}\hskip 1.5pt}

\newmuskip\pFqskip
\mathchardef\pFcomma=\mathcode`, 

\numberwithin{equation}{section}

\begin{document}

\title{\bf 

A note on $\mathfrak{su}(2)$ models and the biorthogonality of generating functions of Krawtchouk polynomials
}
\author{

Luc Vinet\textsuperscript{$1$}\footnote{
E-mail: vinet@CRM.UMontreal.CA}~, 
Alexei Zhedanov\textsuperscript{$2$}\footnote{
E-mail: zhedanov@yahoo.com} \\[.5em]
\textsuperscript{$1$}\small~Centre de Recherches Math\'ematiques, 
Universit\'e de Montr\'eal, \\
\small~P.O. Box 6128, Centre-ville Station, Montr\'eal (Qu\'ebec), 
H3C 3J7, Canada.\\[.9em]

\textsuperscript{$2$}\small~School of Mathematics, 
Renmin University of China, Beijing, 100872, China
}
\date{\today}
\maketitle

\hrule
\begin{abstract}\noindent 
Eigenvalue problems on irreducible $\mathfrak{su}(2)$ modules and their adjoints are considered in the Bargmann, Barut-Girardello and finite difference models. The biorthogonality relations that arise between the corresponding generating functions of the Krawtchouk polynomials are sorted out. A link with Pad\' e approximation is made.
\end{abstract}

\bigskip

\hrule

\section{Introduction}

As a rule, the functions that respectively solve (generalized) eigenvalue problems and their transpose will be biorthogonal (see for example \cite{zhedanov1999biorthogonal}). Of course if the operators involved are self-adjoint the solutions associated to different eigenvalues are simply orthogonal.  In the investigation of the algebraic description of certain families of biorthogonal functions \cite{vinet2021unified}, \cite{vinet2021algebraic}, we were led to observations pertaining to the generating functions of the Krawtchouk polynomials \cite{koekoek2010hypergeometric} that prompted this note.

The interpretation of the Krawtchouk polynomials as matrix elements of $SU(2)$ representations \cite{koornwinder1982krawtchouk}, \cite{zhedanov1999biorthogonal}, \cite{floreanini1993quantum}, is possibly the simplest intance of connection between groups, algebras and special functions. It hence offers a nice framework to illustrate results that generalize to more involved situations. Henceforth we shall consider the irreducible representations of $\mathfrak{su}(2)$ and their adjoints. After recalling how the Krawtchouk polynomials appear in this picture, we shall consider three familiar models: the finite difference realization \cite{nikiforov1991classical}, the Bargmann model \cite{gelfand2018representations},\cite{naimark2014linear}, \cite{miller1968lie} and the Barut--Girardello one \cite{barut1971new}, \cite{de2019Bargmann}. In each case we shall describe how biorthogonality between a natural eigenvalue problem on irreducible modules and its transpose arises. For the latter two realizations this will amount to biorthogonality relations between generating functions; furthermore a connection with Pad\' e approximation will also be underscored.

\section{Krawtchouk polynomials}

Drawing from \cite{koekoek2010hypergeometric}, we shall record for convenience some properties of the Krawtchouk polynomials $K_n(k;p, N)$. These polynomials are defined in terms of the Gauss hypergeometric series by
\begin{equation}
    K_n(k;p, N) = {_2}F_1 \left( {-n, -k \atop -N } ;\frac{1}{p}\right), \quad n=1, \dots, N. \label{def}   
\end{equation}
For the points we wish to make, nothing essential will be lost by considering the symmetric case where $p=\frac{1}{2}$ and this will make the formulas lighter. The polynomials satisfy then the orthogonality relation:
\begin{equation}
    \sum_{k=0}^N \binom{N}{k} K_m(k;\frac{1}{2}, N)K_n(k;\frac{1}{2}, N) = 2^N \frac{(-1)^n n!}{(-N)_n} \delta_{mn}.\label{orth}
\end{equation}
They also obey the three-term recurrence relation:
\begin{equation}
    (N-2k) K_n(k;\frac{1}{2}, N)=(N-n) K_{n+1}(k;\frac{1}{2}, N) + n K_{n-1}(k;\frac{1}{2}, N), \label{rec1}
\end{equation}
or with
\begin{equation}
    p_n(k) = \big{(}\frac{1}{2}\big{)}^n (-N)_n \;K_n(k;\frac{1}{2}, N), \label{pn}
\end{equation}
the normalized one:
\begin{equation}
    k p_n(k) = p_{n+1}(k) + \frac{1}{2} N p_n(k) + \frac{1}{4} n(N + 1 -n) p_{n-1}(k). \label{rec2}
\end{equation}

\section{The irreducible representations of $\mathfrak{su}(2)$ and their adjoints}

The Lie algebra $\mathfrak{su}(2)$ has $\{ J_0, J_{\pm}\}$ as generators with relations
\begin{equation}
    [J_0, J_{\pm}]= \pm J_{\pm}, \qquad [J_+, J_-]= 2J_0. \label{cr}
\end{equation}
Its Casimir element is
\begin{equation}
    J^2 = J_0^2 - J_0 + J_+J_-.\label{Cas}
\end{equation}
The irreducible representations of this algebra are well known \cite{gelfand2018representations} to be of dimension $N+1$ with $N$ a positive integer. They are characterized by the fact that in such a representation 
\begin{equation}
    J^2=\frac{N}{2}\big{(}\frac{N}{2} + 1\big{)} \label{Casval}
\end{equation}
on the corresponding modules. A basis for these representations is given by the eigenvectors $\{|n\rangle\}$ of $J_0$ with eigenvalues $-\frac{N}{2}, -\frac{N}{2} + 1, \dots, \frac{N}{2}$. The actions of the generators on these basis vectors are:
\begin{align}
    J_0 |n\rangle =& \;(n - \frac{N}{2})|n\rangle \label{irrep1} \\
    J_+ |n\rangle =& \;(n-N) |n+1\rangle\label{irrep2}\\
    J_- |n\rangle =& - n \;|n-1\rangle. \label{irrep3}
\end{align}

The normalized dual basis $\{\widetilde{|n\rangle}\}$ such that $\widetilde{\langle m|}n \rangle = \delta_{mn}$ consists in the eigenvectors of the transpose $J_0^T$ of $J_0$ with eigenvalue $n-\frac{N}{2}$. This simply follows from 
\begin{equation}
    \widetilde{\langle m|}J_0|n\rangle = (n-\frac{N}{2})\widetilde{\langle m|}n\rangle = (\widetilde{\langle m|}J_0^T)|n\rangle = (m-\frac{N}{2})\widetilde{\langle m|}n\rangle.
\end{equation}
With $\widetilde{\langle m|}X|n\rangle = \langle n|X^T\widetilde{|m\rangle}$, the action of the transposed operators is seen to be
\begin{align}
    J_0^T \widetilde{|n\rangle} =& \;(n - \frac{N}{2})\widetilde{|n\rangle} \label{irreptilde1}\\
    J_+^T \widetilde{|n\rangle} =& \;(n-1-N) \widetilde{|n-1\rangle} \label{irreptilde2} \\
    J_-^T \widetilde{|n\rangle} =& -(n+1) \;\widetilde{|n+1\rangle}. \label{irreptilde3}
\end{align}
The transposition manifestly turns $J_+^T$ into a lowering operator and reciprocally $J_-^T$ into a raising one. In order to have a perfect correspondence between the initial representation given by \eqref{irrep1}, \eqref{irrep2}, \eqref{irrep3} and its transpose we shall renormalize the vectors $\{\widetilde{|n\rangle}\}$ according to
\begin{equation}
|n\rangle^* = \frac{n!(N-n)!}{N!} \widetilde{|n\rangle.} \label{normal}
\end{equation}
This modifies the actions \eqref{irreptilde1}, \eqref{irreptilde2} and \eqref{irreptilde3} into
\begin{align}
    J_0^T |n\rangle^* =& \;(n - \frac{N}{2})|n\rangle^* \label{irrep*1}\\
    J_+^T |n\rangle^* =& \;-n |n-1\rangle^* \label{irrep*2}\\
    J_-^T |n\rangle^* =& \;(n-N) |n+1\rangle^*. \label{irrep*3}
\end{align}
The match between the representation on $span \{|n\rangle\}$ and the one on $span \{|n\rangle^*\}$ is then simply obtained by taking $\Tilde{J_0} = J_0^T$, $\Tilde{J_+}= J_-^T$ and $ \Tilde{J_-}= J_+^T$ as expected from the effect of the transposition on the commutation relations. Observe that the normalization \eqref{normal} is singular for $n<0$ and $n>N$ with the effect of truncating the actions of $J_+^T$ and $J_-^T$.

\section{$\mathfrak{su}(2)$ and the Krawtchouk polynomials}

A simple way to establish the connection that the symmetric Krawtchouk polynomials have with $\mathfrak{su}(2)$ is to consider the generator $X=\frac{1}{2}(J_+ + J_-)$ and to examine its diagonalization on the $span \{|n\rangle\}$ . It is immediate that $X$ can be obtained from $J_0$ by a specific adjoint action of the group $SU(2)$ on its algebra. Hence $X$ will have the same spectrum as $J_0$. We thus posit the eigenvalue problem\footnote{To obtain the Krawtchouk polynomials with an arbitrary $p$ would simply require adding the operator $J_0$ to this $X$}
\begin{equation}
  X|\lambda_k\rangle =   \frac{1}{2}(J_+ + J_-)|\lambda_k\rangle = (k - \frac{N}{2}) |\lambda_k\rangle. \label{EP}
\end{equation}
Now expand $|\lambda_k\rangle$ over the basis $|n\rangle$:
\begin{equation}
    |\lambda_k\rangle = \sum_{n=0}^N C_n(k) |n\rangle. \label{exp}
\end{equation}
Using the actions \eqref{irrep2} \eqref{irrep3}, it is readily seen that the eigenvalue equation \eqref{EP} leads to the following recurrence relation for the coefficients $C_n(k)$:
\begin{equation}
    (N-2k) C_n(k)= (n+1) C_{n+1}(k) + (N+1-k)C_{n-1}(k). \label{recCn}
\end{equation}
Setting 
\begin{equation}
    C_n(k)=(-1)^n\frac{2^n}{n!}p_n(k) 
\end{equation}
brings the identification with \eqref{rec2} and in view of \eqref{pn}, we have
\begin{equation}
    C_n(k)= \binom{N}{n} K_n(k;\frac{1}{2}, N) \label{coeff}
\end{equation}
where we used
\begin{equation}
    (-N)_n=(-1)^n \frac{N!}{(N-n)!}.
\end{equation}
This is all quite familiar. Let us now consider the adjoint eigenvalue problem
\begin{equation}
    X^T |\lambda_k\rangle^* = \frac{1}{2}(J_+^T + J_-^T) |\lambda_k\rangle^* = (k- \frac{N}{2}) |\lambda_k\rangle^*. \label{adjoint}
\end{equation}
Let 
\begin{equation}
    |\lambda_k\rangle^* = \sum_{n=0}^N C^*_n(k) |n\rangle^*. \label{lambda*}
\end{equation}
Since the actions of $J_+^T$ and of $J_-^T$ in the basis $\{|n\rangle^*\}$ coincide respectively with those of $J_-$ and of $J_+$ in the basis $\{|n\rangle\}$, the eigenvalue equation \eqref{adjoint} will yield for the coefficients $C^*_n(k)$ the same recurrence relation as the one, \eqref{recCn}, satisfied by the coefficients $C_n(k)$ introduced before. Hence, given \eqref{normal},
\begin{equation}
    |\lambda_k\rangle^* = \sum_{n=0}^N  \binom{N}{n} K_n(k;\frac{1}{2},N) |n\rangle^* = \sum_{n=0}^N  K_n(k;\frac{1}{2},N) \widetilde{|n\rangle}. \label{lam} 
\end{equation}
\textit{The Krawtchouk polynomials thus appear as the overlaps between the eigenstates of $J^T_0$ and those of $X$ according to \eqref{coeff} and as well, in the overlaps between the eigenstates of $J_0$ and $X^T$ as per \eqref{lambda*}.}

From general linear algebra theory, as solutions of adjoint eigenvalue problems, the vectors $|\lambda_k\rangle$ and $|\lambda_l\rangle^*$ associated to different eigenvalues should be orthogonal. This is readily checked using the orthogonality of the Krawtchouk polynomials and their duality property namely,  $K_n(k; \frac{1}{2}, N) = K_k(n; \frac{1}{2}, N)$ which follows from their definition. Indeed recalling that $\widetilde{\langle m|} n\rangle = \delta_{mn}$, we have
\begin{align}
    ^*\langle \lambda_k | \lambda_l \rangle &=
    \sum_{m,n=0}^N \binom{N}{n}K_m(k; \frac{1}{2}, N)K_n(l; \frac{1}{2}, N) \widetilde{\langle m |}n\rangle\\
    &= \sum_{n=0}^N \binom{N}{n}K_n(k; \frac{1}{2}, N)K_n(l; \frac{1}{2}, N)\\
    &= \sum_{n=0}^N \binom{N}{n}K_k(n; \frac{1}{2}, N)K_l(n; \frac{1}{2}, N) \propto \delta_{kl}. \label{biorth}
    \end{align}
    
The remainder of this note indicates how this is realized in common models of $\mathfrak{su}(2)$ and involves various generating functions of the Krawtchouk polynomials.  

\section{The finite difference model}

Let $T_{\pm}$ be the shift operators acting as follows on functions of $k$:
\begin{equation}
    T_{\pm} f(k) = f(k\pm 1).
\end{equation}
It is readily checked that the assignment
\begin{align}
    J^{(\Delta)}_0 &= s- \frac{N}{2}\\
    J^{(\Delta)}_+ &= -(N-s+1)T_-\\
    J^{(\Delta)}_- &= -(s+1) T_+,
    \end{align}
satisfy the commutation relations \eqref{cr}. In this model, the eigenstates of $J_0$ are realized by 
\begin{equation}
    \widetilde{\langle s|}n\rangle = \delta_{ns}.
\end{equation}
Indeed 
\begin{equation}
    \widetilde{\langle s|}J_0 |n\rangle = J^{(\Delta)}_0 \widetilde{\langle s|}n\rangle = (s- \frac{N}{2}) \delta_{sn} = (n- \frac{N}{2}) \delta_{sn.}
\end{equation}
Since 
\begin{equation}
   \widetilde{\langle s|}\lambda_k\rangle = \sum_{n=0}^N  C_n(k) \widetilde{\langle s|}n\rangle = C_s(k),
\end{equation}
the states $|\lambda _k\rangle$ are represented by the functions $\lambda_k(s)= \widetilde{\langle s|}\lambda_k\rangle$ which from \eqref{coeff} are seen to be the Krawtchouk polynomials themselves. This is verified by observing that
\begin{equation}
   \widetilde{\langle s|}\frac{1}{2}(J_+ + J_-)|\lambda_k\rangle  =  \frac{1}{2}(J^{(\Delta)}_+ + J^{(\Delta)}_-)\lambda_k (s) = \frac{1}{2} \big{[}  - (N-s+1)\lambda_k (s-1) -(s+1)\lambda_k (s+1)\big{]}.
\end{equation}
Hence the eigenvalue equation
\begin{equation}
    \frac{1}{2}(J^{(\Delta)}_+ + J^{(\Delta)}_-)\lambda_k (s) = (k - \frac{N}{2}) \lambda_k(s) \label{diff1}
\end{equation}
is identified with the relation \eqref{recCn} (with index and variable interchanged) to confirm that 
\begin{equation}
    \lambda_k (s) = \binom{N}{s} K_k(s; \frac{1}{2}, N).
\end{equation}
The transposed operators $J^{(\Delta)}_0, J^{(\Delta)}_{\pm}$ are
\begin{align}
    {J^{(\Delta)}_0}^T &= s- \frac{N}{2}\\
    {J^{(\Delta)}_+}^T &= -(N-s+)T_+\\
    {J^{(\Delta)}_-}^T &= -s T_-.
    \end{align}
    with $T^T_{\pm} = T_{\mp}.$
The adjoint eigenvalue problem for $\lambda^*_k (s) = \langle s |\lambda_k\rangle^*$,
\begin{equation}
    \frac{1}{2}({J^{(\Delta)}_+}^T + {J^{(\Delta)}_-}^T) \lambda^*_k (s) = (k - \frac{N}{2}) \lambda^*_k (s)
\end{equation}
translates to
\begin{equation}
    (N-2s) \lambda^*_k(s+1) + s \lambda^*_k(s-1) = (N-2k) \lambda^*_k(s) \label{diff2}
\end{equation}
which shows \cite{nikiforov1991classical}, as expected, when comparing with the recurrence relation (equivalently the difference equation because of the duality symmetry) \eqref{rec1} that
\begin{equation}
    \lambda^*_k(s)=K_k(s; \frac{1}{2}, N).
\end{equation}
We naturally observe the correspondance with the normalization relation \eqref{normal}.
The biorthogonality between the solution $\lambda_k(s)$ of \eqref{diff1} and the solutions $\lambda^*_k(s)$ of the adjoint problem \eqref{diff2}, namely,
\begin{equation}
    \sum_{s=0}^N \lambda^*_k(s) \lambda_l(s) \propto \delta_{kl} 
\end{equation}
therefore follows from the orthogonality of the Krawtchouk polynomials.
\section{The Bargmann model}

We shall discuss next two differential realizations. The so-called Bargmann model has the $\mathfrak{su}(2)$ generators represented by the following differential operators \cite{naimark2014linear} acting on functions of the variable $z$:
\begin{align}
    J^{(B)}_0 = &\;z\partial_z - \frac{N}{2}\\
    J^{(B)}_+ = &\;z^2\partial_z - Nz\\
    J^{(B)}_- = & -\partial_z.
    \end{align}
That the $\mathfrak{su}(2)$ commutation relations \eqref{cr} are satisfied by these operators is easily checked. The basis states $\{|n\rangle\}$ are modelled by the monomials
\begin{equation}
    \widetilde{\langle z |} n\rangle = z^n , \qquad n =0,\dots, N\label{zn}
\end{equation}
since indeed
\begin{equation}
    \widetilde{\langle z |}J_0 | n\rangle = J^{(B)}_0\widetilde{\langle z |} n\rangle = (z\partial_z - \frac{N}{2}) z^n = (n - \frac{N}{2})\widetilde{\langle z |} n\rangle
\end{equation}
and
\begin{equation}
\widetilde{\langle z |}J_+ | N\rangle = J^{(B)}_+\widetilde{\langle z |} N\rangle = 0 \quad \text{and} \quad  \widetilde{\langle z |}J_- | 0\rangle = J^{(B)}_-\widetilde{\langle z |} 0\rangle=0.
\end{equation}
The eigenfunctions $\lambda_k (z) = \widetilde{\langle z |} \lambda_k\rangle$ will be represented by
\begin{equation}
   \lambda_k(z) = \sum_{n=0}^N \binom{N}{n} K_n(k; \frac{1}{2}, N) z^n \label{lamz}
\end{equation}
in view of \eqref{exp}, \eqref{coeff} and \eqref{zn}. They will satisfy 
\begin{equation}
    \widetilde{\langle z |}\frac{1}{2}(J_+ + J_-) | \lambda_k \rangle = \frac{1}{2}(J^{(B)}_+ + J^{(B)}_- )\lambda_k(z) = (k - \frac{N}{2}) \lambda_k (z),
\end{equation}
that is 
\begin{equation}
    \big{[} (z^2 - 1)\partial_z - Nz + (N - 2k)\big{]} \lambda_k (z) = 0.
\end{equation}
This differential equation is readily integrated to obtain a solution on span $\{z^n, n= 0,\dots, N\}$. Indeed one finds 
\begin{equation}
    \lambda_k(z) = \alpha (1-z)^k (1+z)^{N-k},
\end{equation}
a polynomial of degree $N$, with $\alpha$ an integration constant which will be set equal to 1 to ensure the match with \eqref{lamz}. This yields a simple and well known derivation of the generating formula:
\begin{equation}
    (1-z)^k (1 + z)^{N-k} = \sum_{n=0}^N \binom{N}{n} K_n(k; \frac{1}{2}, N) z^n. \label{gen1}
\end{equation}

Consider now the Lagrange adjoints which read:
\begin{align}
{J^{(B)}_0} ^T =& \;- z\partial_z - \frac{N}{2} - 1\\
{J^{(B)}_+} ^T =& \;- z^2\partial_z - (N+2)z\\
{J^{(B)}_-} ^T =& \; \; \partial_z.
\end{align}
 The eigenfunctions $\langle z \widetilde{|n\rangle}$ of ${J^{(B)}_0} ^T$ are
    \begin{equation}
        \langle z \widetilde{|n\rangle} = z^{-n - 1}, \qquad n=0,\dots, N,
    \end{equation}
    as is seen from
  \begin{equation}
      \langle z |J_0^T \widetilde{| n \rangle} = {J^{(B)}_0} ^T \langle z \widetilde{|n\rangle} = - (z\partial_z + \frac{N}{2} + 1) z^{-n-1} = (n- \frac{N}{2})\langle z \widetilde{|n\rangle}.
  \end{equation}  
  We observe that the eigenfunctions $\widetilde{\langle z |} n\rangle$ and $\langle z \widetilde{|n\rangle}$ of $J^{(B)}_0 $ and of ${J^{(B)}_0 }^T$ are orthogonal with respect to the scalar product provided by integration in the complex plane along a contour $\Gamma$ encircling the origin:
  \begin{equation}
    \frac{1}{2\pi i}\oint_\Gamma dz  \widetilde{\langle m |} z \rangle \widetilde{\langle z |}n \rangle = \frac{1}{2\pi i}\oint_\Gamma dz \;z^{-1-m+n} = \delta_{mn}.
  \end{equation} \label{sp}
Formally it is verified that when acting on the functions
\begin{equation}
    \langle z| n\rangle^* = \frac{n!(N-n)!}{N!} z^{-n-1},
\end{equation}
the operators ${J^{(B)}_0} ^T$, ${J^{(B)}_{\pm}} ^T$ imitate the representation given by \eqref{irrep*1}, \eqref{irrep*2}, \eqref{irrep*3}. For instance
\begin{align}
    {J^{(B)}_-} ^T \langle z | n\rangle^* &= -(n+1) \binom{N}{n}^{-1} z^{-n-2}\\
    & = \;-(n+1) \binom{N}{n}^{-1} \binom{N}{n+1} \langle z | n + 1\rangle^*\\
    & = \:(n-N) \langle z | n + 1\rangle^*
\end{align}
in agreement with \eqref{irrep*3}. This seemingly shows that $ {J^{(B)}_-} ^T \langle z | N\rangle^* = 0$. Note however that $\langle z | N\rangle^*= z^{-N-1}$ and hence that $ {J^{(B)}_-} ^T z^{-N-1} = -(N+1)z^{-N-2} \neq 0$. The truncation for $n > N$ must hence be imposed by hand. Giving a priori the function $\langle z | N+1\rangle^*$ an infinite normalization is signalling this restriction. Similarly, the action $ {J^{(B)}_+} ^T \langle z | 0\rangle^* = 0$ must also be declared irrespective of the action of the differential operator. 
It is with this understanding that the eigenvalue problem for ${X^{(B)}}^T = \frac{1}{2}({J^{(B)}_+} ^T  + {J^{(B)}_-} ^T)$ should be set in the Bargmann model. This 
is how the function 
\begin{equation}
    \lambda^*_k(z) = \langle z |\lambda_k\rangle^* = \sum_{n=0}^N K_n(k ; \frac{1}{2} , N) z ^{-n-1} \label{lamad}
\end{equation}
obtained from \eqref{lam} can be viewed as an eigenfunction of 
\begin{equation}
    {X^{(B)}}^T = \frac{1}{2}({J^{(B)}_+} ^T  + {J^{(B)}_-} ^T) = \frac{1}{2} \big{[}(1-z^2)\partial_z - (N+2)z\big{]}
\end{equation}
on $span\{\binom{n}{N}^{-1} z^{-n-1}, n=0,\dots,N\}$ with eigenvalue $(k - \frac{N}{2})$. It should thus be stressed that because the truncations mentioned above must be imposed, the functions $\lambda^*_k(z)$ cannot be obtained by solving freely the differential equation that ${X^{(B)}}^T\lambda^*_k(z) = (k - \frac{N}{2}\lambda^*_k(z)$ would appear to entail.

This can be rephrased as follows. Let $\Pi$ be the projector from the space of Laurent series in $z$ to $span \{z^{-1-n}, n=0,\dots,N\}$. The generators $ J_{\bullet}^T$ of $\mathfrak{su}(2)$ in this model should really be represented by $\Pi({J^{(B)}_{\bullet}}^T) \Pi$. The solution of the equation
\begin{equation}
    2 \Pi\  ({X^{(B)}}^T) \Pi \;\lambda^*_k(z) = \Pi \big{(}\big{[}(1-z^2)\partial_z - (N+2)z\big{]}\big{)}\Pi \;\lambda^*_k(z) = (2k -N) \lambda^*_k(z)
\end{equation}
 can thus be obtained by taking an ansatz of the form
\begin{equation}
    \lambda^*_k(z) = \sum_{n=0}^N a^{(k)}_n z^{-1-n}.
\end{equation}
Substituting in the equation gives
\begin{equation}
    \Pi \sum_{n=0}^N  a^{(k)}_n \big{[} -(n+1)z^{-2-n} + (n-N-1)z^{-n} + (N-2k) z^{-1-n}\big{]} = 0,
\end{equation}
which since $\Pi (z^{0}) = 0$, yields for the coefficients $a^{(k)}_n$ the recurrence relation
\begin{equation}
    -n a^{(k)}_{n-1} + (n-N) a^{(k)}_{n+1} + (N - 2k) a^{(k)}_n = 0.
\end{equation}
This relation is immediately identified with \eqref{rec1} to give $a^{(k)}_n=K_n(k ; \frac{1}{2} , N)$ in conformity with \eqref{lamad}.

The functions $\lambda_k(z)$ and $\lambda^*_l (z)$ prove biorthogonal under the scalar product \eqref{sp}. Indeed,
\begin{align}
    \frac{1}{2\pi i}\oint_\Gamma dz \lambda_k(z)\lambda^*_l (z) &=  \frac{1}{2\pi i}\oint_\Gamma dz \sum_{n,m=0}^N \binom{N}{m}K_n(k ; \frac{1}{2} , N) K_m(l ; \frac{1}{2} , N)z^{m-n-1}\\
    &= \sum_{n=0}^N \binom{N}{n} K_n(k ; \frac{1}{2} , N) K_n(l ; \frac{1}{2} , N) \propto \delta_{kl}.
\end{align}

Let us mention that $\lambda^*_k(z)$ can be given an alternative expression as a truncated series using the generating function  (eq. (1.10.13)) given in \cite{koekoek2010hypergeometric}.
One has that
\begin{equation}
    \lambda^*_k(z) = \bigg{[} \bigg{(}\frac{1}{z-1}\bigg{)}\; {_2}F_1 \left( {1, -k \atop -N } ;\frac{2}{1-z}\right)\bigg{]}_N\label{form}
\end{equation}
where the subscript $N$ means : truncate the power series in $\frac{1}{z}$ that will start with $z^{-1}$ after N terms.
It thus follows that this truncated series is orthogonal (for different eigenvalues) to the product $(1-z)^l (1+ z)^{N-l}. $ Formula \eqref{form} can be checked directly, 
it would however be satisfactory to have an algebraic derivation relying on $\mathfrak{su}(2)$ representation theory. (See \cite{groenevelt2001meixner}, \cite{baeder2015power} in this connection.)

\section{The Barut--Girardello model}

The last model we shall consider is obtained by a Laplace transform from the Bargmann one \cite{de2019Bargmann} and has the names of Barut and Girardello apposed to it. In this section, we shall first look at how the eigenvalue problem for $(J_+ + J_-)$ is realized and discuss how it leads to a generating function involving the confluent hypergeometric series. Second, we shall describe a connection with Pad\' e approximation. Third we shall examine the adjoint eigenvalue problem to obtain the biorthogonal partner of the solutions to the initial problem; we shall indicate how the generating function stemming from the adjoint problem  is equivalent to the previous one owing to the persymmetric properties \cite{genest2017persymmetric} of the Krawtchouk polynomials. 

In the Barut--Girardello model, the generators are again realized as differential operators acting on the variable z; they are:
\begin{align}
    J^{(BG)}_0 = &\;z\partial_z - \frac{N}{2}\\
    J^{(BG)}_+ = &\;z\\
    J^{(BG)}_- = & -z\partial^2_z + N\partial_z.
    \end{align}
    
The elements of the monomial basis $\{ \widetilde{\langle z |} n\rangle = z^n, n=0, \dots, N\}$ need to be given a different normalization for the actions
\eqref{irrep1}, \eqref{irrep2} and \eqref{irrep3} to be reproduced by the operators $J^{(BG)}_0$ and $J^{(BG)}_{\pm}$. Denoting these new basis vectors associated to the Barut--Girardello model by $\widehat{| n\rangle}$, we shall have
\begin{equation}
    \widetilde{\langle z } \widehat{| n \rangle} = \frac{1}{(-N)_n}  \widetilde{\langle z} | n \rangle = \frac{(-1)^n (N-n)!}{N!} z^n.
\end{equation}
While this normalization will generally yield the desired effect as one can check, here again a truncation must be imposed: one must require that $J^{(BG)}_+ \widetilde{\langle z}\widehat{|N\rangle } = 0$ in keeping with the fact that the eventual $\widehat{|N+1}\rangle $ vector would be ``infinite". Note however that $\widetilde{\langle z} \widehat{| 0 \rangle} =1$ is naturally annihilated by $J^{(BG)}_-$.

The eigenvectors $|\lambda_k \rangle$ of $\frac{1}{2}(J_+ + J_-)$ will be represented by the function
\begin{equation}
    \widehat{\lambda_k} (z)= \sum_{n=0}^N D_n^k \:\widetilde{\langle z } \widehat{| n \rangle}
\end{equation}
verifying
\begin{equation}
    ( J^{(BG)}_+ +  J^{(BG)}_-)\widehat{\lambda_k} (z)= (2k-N)\widehat{\lambda_k} (z). \label{EVPBG}
\end{equation}
Since $J^{(BG)}_+$ and $J^{(BG)}_-$ model the actions \eqref{irrep*2} and \eqref{irrep*3} on the functions $\widetilde{\langle z } \widehat{| n \rangle}=\frac{1}{(-N)_n}z^n$, the expansion coefficients $D_n^k$ will coincide with the coefficients $C_n^k$ given in \eqref{coeff}. It thus follows that  $\widehat{\lambda_k} (z)$ is given by
\begin{equation}
  \widehat{\lambda_k} (z) = \sum_{n=0}^N \frac{(-1)^n}{n!} K_n(k; \frac{1}{2}, N) z^n. \label{lamzhat}
    \end{equation}
Mindful of the truncation that must be enforced, this function should satisfy the differential equation that the eigenvalue problem \eqref{EVPBG} implies and which reads:
\begin{equation}
    \big{[}-z \partial_z^2 + N\partial_z + (z + N - 2k)\big{]} \widehat{\lambda_k} (z) = 0.
\end{equation}
The solutions of this equation are expressible in terms of the confluent hypergeometric function ${_1}F_1$. Remembering that it should be restricted to $span\{1, z, \dots, z^N\}$ one finds
\begin{equation}
    \widehat{\lambda_k} (z) = \Bigg{[} e^{-z}\;{_1}F_1 \left( { -k \atop -N } ; 2z\right)\Bigg{]}_N,
\end{equation}
where the subscript $N$ indicates that the power series in $z$ should be truncated after the term $z^N$. Putting this together with \eqref{lamzhat}, we recover in the special case $p=\frac{1}{2}$ another generating function for the Krawtchouk polynomials ( see eq. (1.10.12) in \cite{koekoek2010hypergeometric}), that is
\begin{equation}
    \Bigg{[} e^{z}\;{_1}F_1 \left( { -k \atop -N } ; -2z\right)\Bigg{]}_N = \sum_{n=0}^N \frac{1}{n!} K_n(k; \frac{1}{2}, N) z^n. \label{gen1F1}
\end{equation}

We here want to point out a connection that the above truncated series have with Pad\' e approximation.  In the theory of confluent hypergeometric function the following formula \cite{bateman1953higher}:
\begin{equation}
   {_1}F_1 \left({ a \atop c}; z \right) = e^z \: {_1}F_1 \left({ c-a \atop c}; -z \right)
\label{tr}  
\end{equation}
is well known and plays an important role in many applications. This Kummer transformation is a limiting case of the Euler transformation of the Gauss hypergeometric function ${_2}F_1(z)$ \cite{bateman1953higher}:
\begin{equation}
      {_2}F_1 \left({ a, b \atop c}; z \right) = (1 - z)^{-b} \: {_2}F_1 \left({ c-a, b \atop c}; \frac{z}{z-1} \right).
\label{euler}  
\end{equation}
These transformations are valid if $c$ is not a negative integer, otherwise the hypergeometric series in \eqref{tr} or \eqref{euler} are not well defined.
One can go around this problem in the following way. Suppose that $c$ is a negative integer, i.e
\begin{equation}
 c=-N, \quad N=1,2,\dots .
\label{c=-N}   
\end{equation}
following the previous notation, we consider
\begin{equation}
\Bigg{[}{_1}F_1 \left({ a \atop -N}; z \right) \Bigg{]}_N =  \sum_{k=0}^N \frac{(a)_k}{k! (-N)_k}z^k
\label{tr_11} 
\end{equation}
that is,  we truncate the hypergeometric summation just before a singularity will appear. Then, formula \eqref{tr} should be replaced by
\begin{equation}
\Bigg{[}{_1}F_1 \left({ a \atop -N}; z \right)\Bigg{]}_N- e^z \: \Bigg{[}{_1}F_1 \left({ -N-a \atop -N}; -z \right)\Bigg{]}_N = O(z^{N+1});
\label{tr_N} 
\end{equation}
in other words, the identity \eqref{tr} remains true if the series are restricted to terms up to $z^{N}$.
The proof of this result is easy and follows the lines of the demonstration given in \cite{vinet2021unified} of the analog for the Euler formula.

Consider the special case of the relation \eqref{tr_N} that occurs when the parameter $a$ is a non-positive integer bigger or equal to $-N$ and let $N= m + n$ with $m$ non-negative
The identity \eqref{tr_N} is then directly seen to imply that the rational function
\begin{equation}
R_{nm}(z) = \frac{{_1}F_1 \left({ -n \atop -n-m}; z \right)}{{_1}F_1 \left({ -m \atop -n-m}; -z \right)} 
\label{R_nm} 
\end{equation}
provides the Pad\'e approximation of the exponential function $e^z$ to within terms $O(z^{n+m})$, that is,
\begin{equation}
R_{nm}(z) -e^z = O(z^{n+m+1}).
\label{R_exp} 
\end{equation}
This gives the well known formula for the Pad\'e table of $e^z$ \cite{baker1996pade}.
We thus see that the generating function of Krawtchouk polynomials that arises from the Barut--Girardello model is related to the Pad\'e approximation of the exponential function. 

We finally come to the transpose of the the eigenvalue problem \eqref{EVPBG}. The Lagrange adjoints of $ J^{(BG)}_0$ and $ J^{(BG)}_{\pm}$ are:
\begin{align}
    {J^{(BG)}_0} ^T = &\;- z\partial_z - \frac{N}{2}- 1, \label{BG1}\\
    {J^{(BG)}_+} ^T= &\;z, \label{BG2}\\
    {J^{(BG)}_-} ^T  = & -z\partial^2_z - (N + 2)\partial_z.\label{BG3}
    \end{align}

The eigenfunctions of ${J^{(BG)}_0} ^T$ are the same as those of ${J^{(B)}_0} ^T$, namely $ \langle z \widetilde{|n\rangle} = z^{-n - 1}$, since these two operators coincide. These eigenfunctions were observed to be orthogonal to  $\widetilde{\langle z|}n\rangle = z^n$ in the last section with the scalar product defined by integration on a contour around the origin.

Another normalization of these basis functions needs to be introduced in order for the operator \eqref{BG1}, \eqref{BG2} and \eqref{BG3} to reproduce in the main the action given in \eqref{irrep*1}, \eqref{irrep*2} and \eqref{irrep*3}. In this case the basis vectors $|n\rangle^*$ should be modeled by
\begin{equation}
   \langle z \widehat{| n\rangle}^* =  (-1)^n n! \langle z \widetilde{|n\rangle} = (-1)^n n! z^{-1-n}.
\end{equation}
Here we need to impose by hand that $ {J^{(BG)}_+} ^T \langle z \widehat{| 0\rangle}^* = 0$ observing that the normalization of the vector $\widehat{| -1\rangle}^*$ if it existed, would be ill-defined. Note however that $ {J^{(BG)}_-} ^T \langle z\widehat{| N\rangle}^* = 0$ is satisfied directly. It is interesting to compare the Bargmann and Barut--Girardello models in this respect. In the former the direct realization respects the domain naturally while the adjoint operators need to be truncated at both ends. For the Barut-Girardello model both the direct and the adjoint actions require the enforcing of \textit{one} truncation.

Consider now the solutions $\widehat{\lambda}^*_k(z)$ of the eigenvalue equation 
\begin{equation}
    \frac{1}{2} ( {J^{(BG)}_+} ^T + {J^{(BG)}_-} ^T)\; \widehat{\lambda}^*_k(z) = (k - \frac{N}{2}) \; \widehat{\lambda}^*_k(z). \label{evplast}
\end{equation}
Owing to the choice of normalization, we know that the expansion coefficients of $\widehat{\lambda}^*_k(z)$ in the basis $\{ \widehat{ | n \rangle}^* \}$ will be the coefficients $C_n^k$ given in \eqref{coeff}. We therefore have
\begin{equation}
  \widehat{\lambda}^*_k(z) = \sum_{n=0}^N \frac{(-1)^n N!}{(N-n)!}\; K_n(k; \frac{1}{2}, N) \;z^{-1-n}. \label{lam**}
\end{equation}
Through steps that we have used repeatedly and that rely on the orthogonality and duality properties of the Krawtchouk polynomials, we observe in this case also that the eigenfunctions $\widehat{\lambda_k} (z)$ of $( J^{(BG)}_+ +  J^{(BG)}_-)$ and $\widehat{\lambda}^*_k(z)$ of the adjoint operator $( {J^{(BG)}_+} ^T + {J^{(BG)}_-} ^T)$ are biorthogonal when the scalar product for which the monomials $\{z^n, n\in \mathbb{Z} \}$ are orthogonal in the complex plane is used:
\begin{align}
    \frac{1}{2\pi i}\oint_\Gamma dz \widehat{\lambda}_k(z) \widehat{\lambda}^*_l (z) &=  \frac{1}{2\pi i}\oint_\Gamma dz \sum_{n,m=0}^N \frac{1}{n!} \frac{N!}{(N-n)!}K_n(k ; \frac{1}{2} , N) K_m(l ; \frac{1}{2} , N)z^{m-n-1}\\
    &= \sum_{n=0}^N \binom{N}{n} K_k(n ; \frac{1}{2} , N) K_l(n ; \frac{1}{2} , N) \propto \delta_{kl}.
\end{align}

Looking at the expressions \eqref{lamz}, \eqref{lamad} and \eqref{lamzhat}, \eqref{lam**}, one observes that the Bargmann and Barut--Girardello models lead to pairs of biorthogonal functions given as generating series of the Krawtchouk polynomials in $z^n$ and $z^{-1-n}$ , $n=0,\dots,N$,  with coefficients that factor the weight term in two natural ways: $\big{[} \binom{N}{n} \cdot 1\big{]}$ and $\big{[} \frac{(-1)^n}{n!} \cdot \frac{1}{(-N)_n}\big{]}$.

The function $\widehat{\lambda}^*_k(z)$ will also verify the differential equation stemming from \eqref{evplast} under the truncation restriction mentioned before. This equation reads:
\begin{equation}
 \big{[}-z \partial_z^2 - (N+2)\partial_z + (z + N - 2k)\big{]} \widehat{\lambda}^*_k(z) = 0. 
\end{equation}
Solving and restricting to $span \{ z^{-1-N}, z^{-N}, \dots, z^{-1}\}$ we find 
\begin{equation}
  \widehat{\lambda}^*_k(z) = (-1)^{N-k} N!\; z^{-1-N} \Bigg{[} e^z\;{_1}F_1 \left( { k-N \atop -N } ; -2z\right)\Bigg{]}_N  \label{genl}
\end{equation}
which will provide another generating function for the Krawtchouk polynomials. The $z$- independent factor in front of \eqref{genl} is introduced to ensure the equality of the two expressions  \eqref{lam**} and \eqref{genl} of $ \widehat{\lambda}^*_k(z)$, it can be identified from the coefficient of $z^{-1-N}$. The subscript $N$ indicates again that the power series should be truncated after $N$ terms.

Upon equating the two expressions given for $ \widehat{\lambda}^*_k(z)$ in \eqref{lam**} and \eqref{genl} and after simple operations including a relabelling of the summation index, one arrives at the identity
\begin{equation}
    (-1)^k  \Bigg{[} e^{-z}\;{_1}F_1 \left( {k-N \atop -N } ; 2z\right)\Bigg{]}_N = \sum_{n=0}^N \frac{1}{n!}\; K_{N-n}(k; \frac{1}{2}, N) \;z^n. \label{genfin}
\end{equation} 
We note that it offers a generating function for the Krawtchouk polynomials with mirror-reflected $(n \leftrightarrow{N-n})$ degrees. 

It can be seen that \eqref{genfin} does not bring a new generating relation and that it is equivalent to \eqref{gen1F1}. This can be attributed to the fact that the Krawtchouk polynomials belong to the persymmetric class \cite{genest2017persymmetric} and that as such the values of $K_n(k; \frac{1}{2}, N)$ and of 
$K_{N-n}(k; \frac{1}{2}, N)$ on the spectral points $k = 0, \dots, N$ are related and specifically verify \cite{borodin2002duality}:
\begin{equation}
    K_{N-n}(k; \frac{1}{2}, N) = (-1)^k K_n(k; \frac{1}{2}, N), \qquad k= 0, \dots, N.\label{mirror}
\end{equation}
This property can be obtained for instance by using the restriction \cite{vinet2021unified} of the transformation formula \eqref{euler} which is required when $c$ is a negative integer together with the explicit expression \eqref{def}. Hence, with the help of \eqref{mirror}, the generating identity \eqref{genfin} can be recast in the form:
\begin{equation}
    \Bigg{[} e^{-z}\;{_1}F_1 \left( {k-N \atop -N } ; 2z\right)\Bigg{]}_N = \sum_{n=0}^N \frac{1}{n!}\; K_n(k; \frac{1}{2}, N) \;z^n. \label{genlast}
\end{equation} 
 This indicates that
 \begin{equation}
     \Bigg{[} e^{-z}\;{_1}F_1 \left( {k-N \atop -N } ; 2z\right)\Bigg{]}_N =  \Bigg{[} e^{z}\;{_1}F_1 \left( { -k \atop -N } ; -2z\right)\Bigg{]}_N ,
 \end{equation}
which is indeed directly implied by \eqref{tr_N} and shows the equivalence of the two generating functions that have arisen in the context of the Barut--Girardello model. One may furthermore take the alternative viewpoint that the mirror-symmetry of the Krawtchouk polynomials also follows from the restricted Kummer transformation formula.
\section{Final remarks}

We have used the simple case of $\mathfrak{su}(2)$ to illustrate issues that relate to the biorthogonality of the solutions of eigenvalue problems and their adjoints on representation spaces. The biorthogonality of the eigenfunctions of a $\mathfrak{su}(2)$ algebra element and of those of the adjoint of this element was observed  to be tantamount to the orthogonality of the Krawtchouk polynomials. In the finite difference model, the adjoint problems correspond to a pair of difference equations solved by the Krawtchouk polynomials themselves. In differential models we noted that attention should be paid to the projections that are required to restrict the action of the generators to the appropriate spaces. This in general precludes solving freely the differential equations that seemingly realize the eigenvalue problems. Normalizations were also seen to play a key role and a connection with Pad\' e approximation was made. The bearing of these considerations on the generating functions of the Krawtchouk polynomials was studied bringing to the fore some biorthogonality properties.

Similar considerations will apply to the Meixner and Charlier polynomials that are respectively related to the $\mathfrak{su}(1,1)$ and the oscillator algebra. In fact, the observations made here will have parallels in the algebraic descriptions of the bispectral polynomials of the Askey scheme \cite{koekoek2010hypergeometric} as well as in the picture in terms of meta-algebras that is being developed  \cite{vinet2021unified} \cite{vinet2021algebraic} for biorthogonal polynomials and rational functions.

\section*{Acknowledgments}
The authors thank H. Cohl for useful correspondence. The work of LV is supported in part by a Discovery Grant from the Natural Sciences and Engineering Research Council (NSERC) of Canada. AZ who is funded by the National Foundation of China (Grant No.11771015) gratefully acknowledges the long-term hospitality of the CRM and the award of a Simons CRM professorship.
\bibliographystyle{unsrt} 
\bibliography{ref_gen.bib}

\end{document}